\documentclass[12pt, reqno]{amsart}
\usepackage{amsmath, amsthm, amscd, amsfonts, amssymb, graphicx, color}

\newtheorem{thm}{Theorem}[section]
\newtheorem{cor}[thm]{Corollary}
\newtheorem{lem}[thm]{Lemma}
\newtheorem{prop}[thm]{Proposition}
\theoremstyle{definition}
\newtheorem{defn}[thm]{Definition}
\newtheorem{exam}[thm]{Example}
\newtheorem{rem}[thm]{Remark}
\numberwithin{equation}{section}

\newcommand{\A}{\mathcal{A}}
\newcommand{\f}{\mathcal{F}}

\newcommand{\R}{\mathcal{R}}

\begin{document}
\title[]{Nonstandard Topology without Logic,\\
Ultrafilters as infinitesimal points in topological space }
\author[M. Akbari Tootkaboni]{M. Akbari Tootkaboni \\
Department of Mathematics, Faculty of science,\\
 Shahed University of Tehran\\
Tehran, Iran\\
Email : akbari@shahed.ac.ir}
\begin{abstract}
Nonstandard analysis is very complex, so finding a simple description of infinitesimal points will be
useful. In this paper, ultrafilters as infinitesimal points in a topological space will be proposed, and some
topological concepts is restated by this tools.

\end{abstract}


\maketitle
\section{\textbf{Introduction}}

Nonstandard Analysis initiated by  Newton and Leibniz,  was
accompanied by logical contradictions which advances in mathematical logic in the twentieth century
could resolve them. Nonstandard methods can give a special insight in topological matters as they are mainly a new way to look at old
things. Using a type-theoretical version of higher-order logic,  A. Robinson in 1960 \cite{Robinson}, introduced the notion of an enlargement which is a
main tool in nonstandard topology. An enlargement is a certain kind of nonstandard
model satisfying a sort of saturation property which is closely related  to  an essential feature of nonstandard methods in
topology, i.e., compactness.

Luxemburg  gave a  simplified version of Robinson's theory in terms of higher-order structures and higher-order languages \cite{Luxemburg}. Edward
Nelson  provided internal set theory (IST) which is an axiomatic basis for a portion
of Robinson's nonstandard analysis  \cite{Nelson} ,\cite{salbany}.

The connections between nonstandard extensions and ultrafilters have been
repeatedly considered in the literature, starting from the seminal paper \cite{Luxemburg} by
Luxemburg (see also \cite{Fenstad}, \cite{Haddad1}, \cite{Haddad2}, \cite{Haddad3} and \cite{Puritz}).
In \cite{Nasso}, M. Di Nasso and M. Forti introduced a notion of topological extension of a given set $X$. The resulting class of
topological spaces includes the Stone-$\check{C}$ech compactification $\beta X$ of the discrete space $X$, as well as all
nonstandard models of $X$ in the sense of nonstandard analysis (when endowed with a natural topology). They
gave a simple characterization of nonstandard extensions in purely
topological terms, and they established connections with special classes of ultrafilters whose existence is
independent of ZFC.

In \cite{Hin-Lead}, the ultrafilter semigroup $(0^+,+)$ of the topological semigroup $T=((0,+\infty),+)$ consists of all nonprincipal ultrafilter
on $T=(0,+\infty)$ converging to the $0$ has been described. According to  \cite{Hin-Lead}, in a topological semigroup $T$, \cite{Akbari-vahed}
has been presented.

On base of \cite{Akbari-vahed}, let $(X,\tau)$ be a topological space and $\beta X_d$ denotes the Stone-$\check{C}$ech compactification of $X_d$, where
$X_d$ is $X$ with discrete topology. The collection of all $p\in\beta X_d$ that converge to $x\in X$ is denoted by
$x^*$, i.e.
\[
x^*=\{p\in \beta X_d:x\in\bigcap_{A\in p}cl_X(A)\}.
\]
 Points of $x^*$ are
infinitely close to $x$, and this help us to obtain some topological properties that has been stated in nonstandard topology. This definition
is explained in Section 3. Also in this section we present some topological properties by
points of $x^*$. In Section 4, we describe separation axioms and study generated topology by this method.

\section{\bf Preliminary}

Let $\Gamma$ be a family of sets that together with $A$ and $B$ contains the intersection $A\cap B$.
By a filter in $\Gamma$ we mean a non-empty subfamily $\f\subseteq \Gamma$ satisfying the following conditions:\\
$(F_1)$ $\emptyset\notin\f$.\\
$(F_2)$ If $A_1,A_2\in\f$, then $A_1\cap A_2\in\f$.\\
$(F_3)$ If $A\in\f$ and $A\subseteq B\in \Gamma$, then $B\in\f$.\\
A filter $\f$ in $\Gamma$ is a maximal filter or an ultrafilter in $\Gamma$, if
for every filter $\A$ in $\Gamma$ that contains $\f$ we have $\A=\f$.

A filter-base in $\Gamma$ is a non-empty family $\mathcal{G}\subseteq \Gamma$ such that
$\emptyset\notin\mathcal{G}$ and
$\mbox{If }A_1,A_2\in\mathcal{G},\mbox{ then there exists }A_3\in\mathcal{G}\mbox{ such that }A_3\subseteq A_1\cap A_2.$
One readily sees that for any filter-base $\mathcal{G}$ in $\Gamma$, the family
\[
\f_{\mathcal{G}}=\{A\in \Gamma:\mbox{ there exists }B\in\mathcal{G}\mbox{ such that }B\subseteq A\}
\]
is a filter in $\Gamma$. For a topological space $(X,\tau)$, a filter in $\tau$ is called open-filter.

A filter $\f$ on $(X,\tau)$ converges to a
point $x\in X$ if $\tau_x\subseteq\f$, where
\[
\tau_x=\{U\subset X:x\in V\subseteq U\mbox{ for some }V\in \tau\}
\]
is the collection of all neighborhoods of $x\in X$.
In this case, the point is called a limit of
the filter $\f$ and we write $x\in lim\f$. A point $x$ is called a cluster point of a filter $\f$ if $x$
belongs to the closure of every member of $\f$. Clearly, $x$ is a cluster point of a filter $\f$ if and only if
every neighborhood of $x$ intersects all members of $\f$. This implies in particular
that every cluster point of an ultrafilter is a limit of this ultrafilter. It is obvious that a subset $A\subseteq X$ is closed in $\tau$ if and only if
a limit of any filter containing $A$ belongs to $A$.

 Let $X_d$ denote $X$ with discrete topology. Now we describe the \mbox{Stone-$\check{C}$ech} compactification $\beta X_d$ of $X_d$.
We take the points of $\beta X_d$ to be the ultrafilters on
 $X_d$,  identifying
the points of  $X_d$ with $\widehat{x}=\{A\subseteq X_d:x\in A\}$ and let $\widehat{A}=\{p\in\beta X_d:A\in p\}$ for $A\subseteq X$. The topology of $\beta X_d$ is defined by
 stating  $\{\widehat{A}:A\subseteq X\}$ as a base for the open sets. Then $cl_{\beta X_d}A=\widehat{A}$ for $A\subseteq X_d$.
 If $\A$ be a filter then $\overline{\A}=\{p\in\beta X_d:\A\subseteq p\}$ is a closed subset of $\beta X_d$. Also, if
 $T$ be a subset of $\beta X_d$, then $\A=\bigcap T$ is a filter and $cl_{\beta X_d}T=\overline{\A}$.
 For more details see \cite{hindbook} and \cite{zelbook}. If $X$ be a completely regular space, then the Stone-$\check{C}$ech compactification of $X$
 has been described by similar way as space of $z$-ultrafilters, see \cite{Gil} and \cite{Akbari}.

 If $X$ and $Y$ are any completely regular spaces, then any continuous function $f:X\rightarrow Y$
 has a unique continuous extension $f^\beta :\beta X\rightarrow \beta Y$.

\begin{lem}
a) Let $f:X_d\rightarrow X_d$ be a function. For each $p\in\beta
X_d$, $$f^\beta(p)=\{A\subseteq X_d: f^{-1}(A)\in p\}.$$ In
particular, if $A\in p$, then $f(A)\in f^\beta(p)$.\\
b) Let $f,g:X_d\rightarrow X_d$ be functions and let $p\in\beta X_d$
satisfy $f^\beta(p)=g^\beta(p)$. Then for each $A\in p$,
$\{x\in X_d:f(x)\in g(A)\}\in p.$\\
c) Let $f:X_d\rightarrow X_d$ be a function. The mapping
$f^\beta:\beta X_d\rightarrow \beta X_d$ has a fixed point if
and only if every finite partition of $X_d$ has a cell $C$ for which
$C\cap f(C)\neq\emptyset.$\\
d) Let $f:X_d\rightarrow X_d$ be injective and have the property
$f(x)\neq x$ for all $x\in X_d$. Then there is a partition
$X_d=A_0\cup A_1\cup A_2$ such that $f(A_i)\cap A_i=\emptyset$ for
$i=0,1,2.$
\end{lem}
{\bf Proof :} For (1), (2) and (3) See \cite{hindbook} Lemma 3.30,3.31 and
3.32. For (4), see \cite{pym}.$\hfill\blacksquare$
\begin{thm}
Let $f:X_d\rightarrow X_d$ be a function. If $f$ has no fixed
points, neither does $f^\beta:\beta X_d\rightarrow \beta X_d$.
\end{thm}
{\bf Proof :} See \cite{hindbook} Theorem 3.34.$\hfill\blacksquare$
\begin{thm}
Let $f:X_d\rightarrow X_d$ be a function and $p\in\beta X_d$. Then
$f^\beta(p)=p$ if and only if $\{x\in X_d:f(x)=x\}\in p$
\end{thm}
{\bf Proof :} See \cite{hindbook} Theorem 3.35.$\hfill\blacksquare$

Now we review the definition of partition regularity and a theorem that connect it with ultrafilters.
   \begin{defn}
  Let $\R$ be a nonempty set of subsets of $X$.  $\R$ is partition regular if and only if
   whenever $\f$ is a finite set of $\mathcal{P}(X)$ ($ \mathcal{P}(X) $ is the set of all subsets of $X$) and $\bigcup\f\in\R$, there exist $A\in\f$ and
    $B\in\R$ such that $B\subseteq A$.
     \end{defn}

\begin{thm}
Let $\R\subseteq \mathcal{P}(X)$  be a nonempty set
and assume $\emptyset\notin\R$.
 Let $$\R^\uparrow=\{B\in \mathcal{P}(X):A\subseteq B\mbox{ for some }A\in\R\}.$$
 Then (a), (b) and (c) are equivalence.\\
 (a) $\R$ is partition regular.\\
 (b) Whenever $\A\subseteq \mathcal{P}(X)$ has the property that every finite nonempty  subfamily of
 $\A$ has an intersection which is in $\R^\uparrow$, there is
  $ \mathcal{U} \in \beta G_d$ such that $\A\subseteq \mathcal{U} \subseteq \R^\uparrow$.\\
 (c) Whenever $A\in\R$, there is  $ \mathcal{U}\in\beta X_d$ such that $A\in \mathcal{U} \subseteq \R^\uparrow$.
 \end{thm}
{\flushleft\bf Proof.} \cite[Theorem 3.11]{hindbook}.
$\hfill\square$
\begin{rem}
a) Let $(X,\tau)$ be a topological space and pick $x\in X$. Then $\R_x=\{A\subseteq X: x\in cl_XA\}$ is partition
regular and $\R_x=\R_x^\uparrow$. Pick $A\in\R_x$ then there exists $p\in\beta X_d$ such that $A\in p\subseteq \R_x$, see Theorem 2.5. It is
obvious that $\tau_x\subseteq \R_x$.\\
b) Let $(X,\tau)$ be a topological space and $F$ be a nonempty subset of $X$. Then
\[
\R_F=\{A\subseteq X:cl_{X}(A)\cap F\neq\emptyset\},
\]
is partition regular and $\R_F=\R_F^\uparrow$.
\end{rem}
Now we quote some prerequisite material from \cite{Ng} for the description of some properties of nonstandard topology, that
has been presented in new way in this paper.\\
Let $(X,\tau)$ be a topological space and $^*X$ be the nonstandard extension of $X$. Given $x\in X$, the monad of $x$ is denoted and defined as
\[
\mu(x)=\bigcap\{^*U:x\in ^*U,U\in\tau\}.
\]
The nearstandard part of $^*X$ is defined as:
\[
ns(^*X)=\bigcup_{x\in X}\mu(x).
\]
The notion of open sets and closed sets can be characterized in terms of monads.
\begin{prop}
Let $A\subseteq X$, then\\
i) $A$ is open if and only if $\mu(x)\subseteq (^*A)$ for each $x\in A$;\\
ii) $A$ is closed if and only if $\mu(x)\cap (^*A)\neq\emptyset$ implies that $x\in A$;\\
iii) For $x\in X$, $x\in cl_X(A)$ if and only if $\mu(x)\cap (^*A)\neq\emptyset$.
\end{prop}
If the monads are disjoint, they form a partition of $ns(^*X)$, producing
equivalence classes on $ns(^*X)$ whose equivalence relation is denoted by $\approx_\tau$,
(infinitely close to each other w.r.t. $\tau$ ), i.e.
\[
\forall x,y\in ns(^*X)(x\approx_\tau y\Leftrightarrow (\exists z\in X(x,y\in\mu(z)\,\,)\,\,\,)\,\,).
\]
The definition of a monad extends to all nearstandard elements in a natural
way: let $x\in ns(^*X)$, we define $\mu(x)$ as $\{y\in ^*X:y\approx_\tau x\}$. So $y\approx_\tau x$ if and only if
$\mu(x)=\mu(y)$, for $x,y\in ns(^*X)$.
\begin{prop}[Monad and separation axioms]
i) A topological space X is Hausdorff if and only if for each distinct
points $x,y\in X$, $\mu(x)\cap\mu(y)=\emptyset$.\\
ii) A topological space X is regular if and only if for each $x\in X$ and every closed subset
$F$ of $X$, $\mu(x)\cap \bigcap_{F\subseteq U\in\tau}(^*U)=\emptyset$.\\
iii) A topological space X is normal if and only if for each disjoint closed subset
$F$ and $K$ of $X$, $\bigcap_{F\subseteq U\in\tau}(^*U)\cap \bigcap_{K\subseteq U\in\tau}(^*U)=\emptyset$.
\end{prop}

\section{\textbf{Nonstandard topology by ultrafilter method}}

Let $(X,\tau)$ be a topological space. For $x\in X$, with respect to $\tau$ on $X$, we
define
$$x^*=\{p\in\beta X_d:x\in\bigcap_{A\in p}cl_X A\}.$$
In fact, $x^*$ is the collection of all ultrafilters converge to $x$. It is obvious that $\widehat{x}\in x^*$. We say $p\in x^*$ is
a near point to $x$. We define $B(X)=\bigcup_{x\in X}x^*$  and $\infty^*=\beta X_d-B(X)$. $p\in B(X)$ is called bounded ultrafilter and
$p\in \infty^*$ is called unbounded ultrafilter. For $F\subseteq X$ define
\[
F^*=\{p\in\beta X_d:p\subseteq \R_F^\uparrow\},\mbox{ see Remark 2.6}.
\]
It is obvious that $x^*\subseteq F^*$, for each $x\in F$.
\begin{lem}
Let $(X,\tau)$ be a Hausdorff  topological space.\\
$a)$ $p\in\infty^*$ if and only if for each $x\in X$ there exists $A\in p$ such that $x\notin cl_X(A)$.\\
$b)$  If $\tau_x\subseteq p$ then $p\in x^*$.\\
$c)$ Let $U\subseteq X$, then $U$ is a neighborhood of $x$ if and only if $U\in p$ for each $p\in x^*$.\\
$d)$ Let $A\subseteq X$. Then $x\in cl_XA$ if and only if $cl_{\beta
X_d}A\cap x^*\neq\emptyset$.\\
$e)$ For each $x\in X$, $\tau_x=\bigcap x^*$ is a filter and $cl_{\beta X_d}(x^*)= \overline{\tau}_x=\bigcap_{U\in \tau_x}cl_{\beta X_d}U$.\\
$f)$ Let $A\subseteq X$, then $x$ is an interior point of $A$ if and only if $x^*\subseteq cl_{\beta X_d}A$.
In particular,
$A$ is open if and only if $x^*\subseteq cl_{\beta X_d}A$ for each $x\in A$.
\end{lem}
{\flushleft\bf Proof.} $a)$ Let $p\in\infty^*$, so $ p \not \in x^* $ for each $ x \in X$.
Hence $ x \not \in \bigcap_{A\in p}cl_X A $. Thus $ x \not \in cl_X A  $ for some $ A \in p $.

Conversely, suppose for each $ x \in X $, there exists $ A \in p $ such that $ x \not \in cl_X A  $. Hence
$ p \not \in B(X) $ and thus $ p \in \infty ^*$.\\
$b)$ Let $ U \in p $ for each $ U \in \tau_x $, thus $ U \cap A \neq \emptyset $ for each $ U \in \tau_x $
and for each $ A \in p $. This implies $ x  \in cl_X A $ for each $ A \in p $. Therefore $ p \in x^* $.

 $c)$ Let $U \in \tau_ x$ and $ p \in x^* $. Since $ U \cap A \neq \emptyset $ for each $ A \in p $, so $ U \in p $.

 Conversely, let $U\in \bigcap_{p\in x^*}p$ and $x\notin int_X(U)$. Then $x\in cl_X(U^c)$ and by Remark 2.6, there
 exists $p\in\beta X_d$ such that $U^c\in p\subseteq \R$. This is a contradiction.\\
$d)$ By Theorem 2.5 and Remark 2.6, there exists $p\in\beta X_d$ such that $A\in p\in x^*$. This implies that
$cl_{\beta X_d}A\cap x^*\neq\emptyset$. Now let $cl_{\beta X_d}A\cap x^*\neq\emptyset$, so there exists $p\in x^*$
such that $p\in cl_{\beta X_d}A$. Thus $A\in p\in x^*$ and so $x\in cl_{X}A$.

 $e)$ and $f)$ are obvious.$\hfill\blacksquare$

For every net $S=\{x_\alpha\}_{\alpha\in I}$ in a topological
space $X$, the family $\f(S)$, consisting of all sets $A\subseteq
X$ with the property that there exists $\alpha_\circ\in I$ such
that $x_\alpha\in A$ whenever $\alpha\geq\alpha_\circ$, is a
filter in the space $X$,( see Theorem 1.6.12 in \cite{Engel}). So
$\{\{x_\alpha:\alpha>\alpha_\circ\}:\alpha_\circ\in I\}$ has the
finite intersection property.
\begin{lem}
Let $X$ be a topological space. Then :\\
(a) Let $\{x_\alpha\}_{\alpha\in I}$  be a net in $X$. If $x_\alpha\rightarrow p$ in $\beta X_d$ for some $p\in x^*$.
 Then $x_\alpha\rightarrow x$ in $(X,\tau)$.\\
 (b) Let $\{x_\alpha\}_{\alpha\in I}$  be a net in $X$ and let $x_\alpha\rightarrow x$ in $X$. Then there exists $p\in x^*$ that is a cluster point of
$\{x_\alpha\}_{\alpha\in I}$ in $\beta X_d$.\\
(c) Let $A\subseteq X$ be closed. Then $A$ is compact if and only if $cl_{\beta X_d}A\cap\infty^*=\emptyset$.
\end{lem}
{\bf Proof :} $(a)$ Let $U$ be an open neighborhood of $x$. Then $U\in p$,(Lemma 3.1), and so
there exists $\beta \in I$ such that $x_\alpha\in U$ for each $\alpha>\beta$. This implies $x_\alpha\rightarrow x$ in $(X,\tau)$.\\
$(b)$ $\{\{x_\gamma:\gamma>\beta\}:\beta\in I\}$ has the finite intersection property
then there exists $p\in\beta X_d$ that $\{x_\gamma:\gamma>\beta\}\in p$ when $\beta\in I$. Since $x_\alpha\rightarrow x$ in $X$. Thus for each
open neighborhood $U$ of $x$, there exists $\beta \in I$ such that $x_\gamma\in U$ when $\gamma >\beta$ and so $U\in p$. By Lemma 3.1 ii), implies
$p\in x^*$.  It is obvious $p$ is cluster point of $\{x_\alpha\}_{\alpha\in I}$ in $\beta X_d$.

$c)$ Let $A\subseteq X$ is compact and $cl_{\beta X_d}A\cap\infty^*\neq\emptyset$, so
 there exists a net $\{x_\alpha\}_{\alpha\in I}$ in $A$ such that $x_\alpha\rightarrow p$
 for some $p\in cl_{\beta X_d}A\cap\infty^*$. Since $A$ is compact so there is a subnet
 $\{x_\beta\}$ such that $x_\beta\rightarrow x\in A$ and so by $(b)$, there is $q\in x^*$
 such that $x_\beta\rightarrow q$. This implies $p=q$, and hence by Lemma 3.1 we have a contradiction.

 Conversely, Let $cl_{\beta X_d}A\cap\infty^*=\emptyset$ and $A\subseteq X$ is not compact. Hence there is a
  net $\{x_\alpha\}_{\alpha\in I}$ in $A$ such that any subnet of $\{x_\alpha\}_{\alpha\in I}$ is divergent in $A$.
  Since $\beta X_d$ is compact, so there is a subnet $\{x_\beta\}$ such that $x_\beta\rightarrow p\in cl_{\beta X_d}A$.
  Also $p\in \infty^*$, because if $p\notin \infty^*$ then $p\in y^*$ for some $y\in X$. So by (a),
  $x_\beta\rightarrow y$ in $X$. This implies that $y\in cl_X A$ and this is a contradiction. Thus
  $p\in cl_{\beta X_d}(A)\cap \infty^*$ and we have a contradiction.
$\hfill\blacksquare$
\begin{thm}[Robinson's Compactness]
Let $(X,\tau)$ be a topological space. Then $A\subseteq X$ is compact
if and only if for every $p\in cl_{\beta X_d}A$ there exists $x\in A$ such that
$p\in x^*$.
\end{thm}
{\bf Proof :} Suppose, to contrary, that A is compact but there is a point $p\in cl_{\beta X_d}A$ such that
$p\notin x^*$ for each $x\in A$. Then every $x\in A$ possesses an open neighborhood $U_x$, with
$U_x\notin p$. Now, since $A$ is compact, the open cover $\{U_x:x\in A\}$ of A has a finite subcover, say
$\{U_1,...,U_n\}$; i.e. $A\subseteq \bigcup_{i=1}^nU_i$. Thus
\[
cl_{\beta X_d}A\subseteq cl_{\beta X_d}(\bigcup_{i=1}^nU_i)=\bigcup_{i=1}^ncl_{\beta X_d}(U_i).
\]
So $p\in cl_{\beta X_d}(U_i)$ for some $i\in\{1,...,n\}$, this implies that $U_i\in p$, and we have a contradiction.

Conversely, let $\{x_\beta\}$ be a net in $A$. Since $\beta X_d$ is compact space,
so there is a subnet $\{x_\beta\}$ such that $x_\beta\rightarrow p\in cl_{\beta X_d}A$.
  Therefore there exists $x\in A$  such that $p\in x^*$. Now by Lemma 3.2(a), $x_\beta\rightarrow x$ in $A$. This complete
  proof.$\hfill\blacksquare$
\begin{thm}
Let $(X,\tau_{_X})$ and $(Y,\tau_{_Y})$ be topological spaces.
Then the following statements are equivalence.\\
a) $f:X\rightarrow Y$ is continuous.\\
b) For each $x\in X$, $f^\beta(x^*)\subseteq (f(x))^*$.\\
c) $f^\beta:B(X)\rightarrow B(Y)$ is well defined and continuous.
\end{thm}
{\bf Proof :} a) implies b). Let $f:X\rightarrow Y$ is continuous, and pick $p\in x^*$. We
must show $f^\beta(p)\in(f(x))^*$. By Lemma 2.1, we have
\[
f^\beta(p)=\{A\subseteq Y:f^{-1}(A)\in p\}.
\]
 Now let $f^{-1}(A)\in p$ for some $A\subseteq Y$, so $x\in cl_X
f^{-1}(A)$. Hence there exists a net $\{x_\alpha\}\subseteq
f^{-1}(A)$ such that $x_\alpha\rightarrow x$ and so
$f(x_\alpha)\rightarrow f(x)$. This implies $f(x)\in cl_Y A$ and hence $f^\beta (p)\in(f(x))^*$.

b) implies a). Let $f^\beta(x^*)\subseteq
(f(x))^*$ for each $x\in X$. Let there exists a net
$\{x_\alpha\}_{\alpha\in I}$ such that $x_\alpha\rightarrow x$ for
some $x\in X$ and $f(x_\alpha)$ is not convergent to $f(x)$. Since
$x_\alpha\rightarrow x$ so for each open neighborhood $U\in
\tau_{_X}$ of $x$ there exists $\beta_U\in I$ such that
$x_\alpha\in U$ for each $\alpha>\beta_U$. Thus
$\A=\{\{x_\alpha:\alpha>\beta_U\}:x\in U\in \tau_{_X}\} $ has the
finite intersection property. Therefore there exists an
ultrafilter $p$ contains $\A$. It is obvious $p\in x^*$,(
because for each $A\in p$ and for each open  neighborhood $U$ of
$x$, we have $A\cap U\neq\emptyset$ so $x\in cl_XA$.)

Since $f(x_\alpha)$ is not convergent to $f(x)$, so there exists a
sub net $\{x_\beta\}$ such that $x_\beta\rightarrow x$ and for
some open neighborhood $U$ of $f(x)$, $f(x_\beta)\notin U$ for
each $\beta$. Since $U$ is an open neighborhood of $f(x)$ so for each
$B\subseteq Y$, if $f(x)\in cl_YB$ then $U\cap B\neq\emptyset$.
This implies $(f(x))^*\subseteq cl_YU$, in particular, $U\in
f^\beta(p)$ and hence $f^{-1}(U)\in p$, (see Lemma 2.1 and 3.1). Now
we have $\{x_\beta:\beta\}\cap f^{-1}(U)\neq\emptyset$, and this
is a contradiction.

a) and b) implies c). It is obvious.

c) implies a). Let $U\in\tau_Y$ then $cl_{\beta Y_d}(\tau_Y)\subseteq cl_{\beta Y_d}(U)$. This implies
$$cl_{\beta X_d}(\tau_X)\subseteq (f^\beta)^{-1}(cl_{\beta Y_d}(U))=cl_{\beta X_d}(f^{-1}(U)).$$
This implies $f^{-1}(U)\in \tau_X$.$\hfill\blacksquare$

\begin{cor}
Let $(X,\tau_{_X})$ and $(Y,\tau_{_Y})$ be topological spaces.
Then $f:X\rightarrow Y$ is continuous if and only if
$f^\beta:B(X)\rightarrow B(Y)$ by $f^\beta(x^*)\subseteq (f(x))^*$ for each $x\in X$ be well define.
\end{cor}
{\bf Proof :} Obvious.$\hfill\blacksquare$
\begin{thm}[Open Mapping]
Let $(X,\tau_X)$ and $(Y,\tau_Y)$ be two topological spaces, let $f:X\rightarrow Y$ be a continuous function and $x^*$ be
a closed subset of $\beta X_d$ for each $x\in X$.
Then $f$ is open if and only if $(f(x))^*\subseteq f^\beta(x^*)$ for each $x\in X$.
\end{thm}
{\bf Proof :} Assume that $f$ is open, and pick $x \in X$. Therefore
$\{f(U):U\in\tau_x\}$ is a neighborhood base at $f(x)$, and thus
\begin{eqnarray*}
f^\beta(x^*)&=&f^\beta(cl_{\beta X_d}(x^*))\\
&=&f^\beta(\bigcap_{U\in\tau_x}cl_{\beta X_d}U)\mbox{ (by Lemma 3.1(e))}\\
&=&\bigcap_{U\in\tau_x}f^\beta(cl_{\beta X_d}U)\\
&=&\bigcap_{U\in\tau_x}(cl_{\beta Y_d}f(U))\mbox{ ($f$ is continuous)}\\
&\supseteq&\bigcap_{V\in\tau_{f(x)}}(cl_{\beta Y_d}V)\\
&=&(f(x))^*.
\end{eqnarray*}
Conversely, assume $(f(x))^*\subseteq f^\beta(x^*)$ for each $x\in X$. To show $f$ is open, pick $U\in\tau_X$ and $y\in f(U)$.
Then $y=f(x)$ for some $x\in U$, and
\[
y^*=(f(x))^*\subseteq f^\beta(x^*)\subseteq f^\beta(cl_{\beta X_d}U)\subseteq cl_{\beta Y_d}f(U).
\]
Thus $f(U)\in p$ for each $p\in y^*$. So $f(U)$ is a neighborhood of $y$, by Lemma 3.1(c). This complete proof.$\hfill\blacksquare$
\begin{thm}
Let $(X,\tau_X)$ and $(Y,\tau_Y)$ be two topological spaces, let $f:X\rightarrow Y$ be a continuous function, and let
$K$ be a compact subset of $X$. Then $f(K)$ is compact.
\end{thm}
{\bf Proof :} Pick $q\in cl_{\beta Y_d}f(K)$, so there exists $p\in cl_{\beta X_d}K$ such that
$f^\beta(p)=q$. Since $K$ is compact, by Theorem 3.3, there exists $x\in K$ such that $p\in x^*$. As $f$ is continuous,
$f^\beta (x^*)\subseteq (f(x))^*$, by Theorem 3.4, and hence $q=f^\beta(p)\in (f(x))^*$. So $f(K)$
is compact, by Theorem 3.3.$\hfill\blacksquare$
\begin{lem}
Let $(X,\tau)$ be a compact topological space. Then $\iota^\beta(x^*)=x$ for each $x\in X$,
where $\iota^\beta:\beta X_d\rightarrow X$ is continuous extension of $\iota:X_d\rightarrow X$.
\end{lem}
{\bf Proof :} Let $p\in x^*$, then $\tau_x\subseteq p$. Let $\{x_\alpha\}_{\alpha\in I}$ be a net in $X_d$ such that
$x_\alpha\rightarrow p\in x^*$ in $\beta X_d$,
then for each $U\in\tau_x$
there exists $\kappa\in I$ such that $x_\gamma\in U$ for each $\gamma>\kappa$. Hence
$\iota^\beta(x_\alpha)\rightarrow x$ and so $\iota^\beta(p)=x$.$\hfill\blacksquare$
\begin{lem}
Let $(X,\tau)$ be a compact topological space and $S$ be a dense subset of $X$. Then :\\
a)  $\iota^\beta(x^*)=x$ for each $x\in S$, where
$\iota^\beta:\beta S_d\rightarrow X$ is continuous extension of $\iota:S_d\rightarrow X$.\\
c) Let $f:S_d\rightarrow X$ be a function. Then $f^\beta(x^*)=f(x)$ for each
$x\in S$, where $f^\beta:\beta S_d\rightarrow X$ is continuous extension of $f$.
\end{lem}
{\bf Proof :} Obvious.$\hfill\blacksquare$

\begin{lem}
Let $(X,\tau_X)$ and $(Y,\tau_Y)$ be completely regular topological spaces. Let $f:X\rightarrow Y$ be a continuous function and
$\overline{f}:\beta X\rightarrow\beta Y$ and $f^\beta: X_d\rightarrow \beta Y_d$ are continuous extension.
 Then :\\
a)  The following diagram is commutes:
\[
\begin{matrix}
\beta X_d &\stackrel{f^\beta}{\longrightarrow} & \beta Y_d\\
\iota^\beta_X{ \downarrow} & \,& \iota^\beta_Y{\downarrow}\\
\beta X &\stackrel{\overline{f}}{\longrightarrow}& \beta Y, \\
\end{matrix}
\]
(i.e. $\overline{f}\circ \iota^\beta_X=\iota^\beta_Y\circ f^\beta$), where
$\iota^\beta_X$ and $\iota^\beta_Y$ are continuous extension of identity function $\iota_X:X_d\rightarrow X$ and
$\iota_Y:Y_d\rightarrow Y$, respectively. \\
b) The following diagram is commutes:
\[
\begin{matrix}
B(X) &\stackrel{f^\beta}{\longrightarrow} & B(Y)\\
\iota^\beta_X{ \downarrow} & \,& \iota^\beta_Y{\downarrow}\\
 X &\stackrel{f}{\longrightarrow}& Y, \\
\end{matrix}
\]
(i.e. $f\circ \iota^\beta_X=\iota^\beta_Y\circ f^\beta$), where
$\iota^\beta_X$ and $\iota^\beta_Y$ are continuous extension of identity function $\iota_X:X_d\rightarrow X$ and
$\iota_Y:Y_d\rightarrow Y$, respectively.
\end{lem}
{\bf Proof :} Obvious.$\hfill\blacksquare$

\begin{thm}
Let $(X,\tau)$ be a topological space, and $f:X\rightarrow X$ be a continuous function. Then
$f$ has a fixed point if and only if there exists $x\in X$ and $p,q\in x^*$ such that $f^\beta(p)=q$.
\end{thm}
{\bf Proof :} If $f$ has a fixed point then there exists $x\in X$ such that $f(x)=x$. Thus $f^\beta(x^*)\subseteq (f(x))^*$ and this complete proof.\\
Conversely, let there exists $x\in X$ and  $p,q\in x^*$ such that $f^\beta(p)=q$. It is obvious that the following diagram is commute.

\[
\begin{matrix}
\beta X_d &\stackrel{f^\beta}{\longrightarrow} & \beta X_d\\
\iota^\beta{ \downarrow} & \,& \iota^\beta{\downarrow}\\
X &\stackrel{f}{\longrightarrow}& X \\
\end{matrix}
\]
By Theorem 3.3 and Lemma 3.6, implies that for $x\in X$
$$f(x)=f(\iota^\beta(q))=\iota^\beta(f^\beta(p))=\iota^\beta(q)=x.\hfill\blacksquare$$

\begin{defn}
Let $\{x_\alpha\}_{\alpha\in I}$ be a net in $(X,\tau)$ and $x\in
X$. For $p\in\beta X_d$  we say $p-lim x_\alpha=x$ if for each
neighborhood $U$ of $x$ and for each $\beta\in I$ implies that $\{x_\alpha:\alpha>\beta\}\cap
U\in p$.
\end{defn}

\begin{lem}
Let $S=\{x_\alpha\}_{\alpha\in I}$ be a net in $X$ and $x\in X$.
Then the statements a) and b) are equivalent and also a) and b) implies c):\\
a) $x_\alpha\rightarrow x$,\\
b) $\bigcap_{\beta\in I}cl_{\beta
X_d}\{x_\alpha:\alpha>\beta\}\subseteq x^*$,\\
c) For each $p\in (\bigcap_{\beta\in I}cl_{\beta
X_d}\{x_\alpha:\alpha>\beta\}$, $p-lim\,\,x_\alpha=x$.
\end{lem}
{\bf Proof :} $a)\rightarrow b)$. Let $p\in \bigcap_{\beta\in I}cl_{\beta
X_d}\{x_\alpha:\alpha>\beta\}$ so there exists a subnet $\{x_\beta\}$ such that
$x_\beta\rightarrow p$ in $\beta X_d$. Since $x_\beta\rightarrow x$ in $X$, by Lemma 3.2 there is
$q\in x^*$ such that a subnet of $\{x_\beta\}$ is convergent to $q$. Therefore $p=q$ and so $p\in x^*$.\\
b)$\rightarrow$ a). Let $\{x_\alpha\}$ is not convergent to $x$.
Therefore there exists a neighborhood $U$ of $x$ such that for some
subnet $\{x_\beta\}_{\beta \in J}$ of $\{x_\alpha\}$ we have
$\{x_\beta:\beta\in J\}\cap U=\emptyset$.

It is obvious $\bigcap_{\gamma\in J}cl_{\beta
X_d}{\{x_\beta:\beta>\gamma\}}\subseteq \bigcap_{\gamma\in I} cl_{\beta
X_d}{\{x_\alpha:\alpha>\gamma\}}$. Now let $p\in cl_{\beta
X_d}{\{x_\beta:\beta>\gamma\}}$. Since  $\tau_x\subseteq p\in x^*$ and
$\{x_\beta:\beta>\gamma\}\in p$
for each $\gamma\in I$, we have a contradiction.\\
a) and b)$\rightarrow$ c). Let $p\in \bigcap_{\beta\in I}cl_{\beta
X_d}\{x_\alpha:\alpha>\beta\}$, then b) implies $p\in x^*$. Thus for
each neighborhood $U$ of $x$ and for each $\beta \in I$, implies $\{x_\alpha:\alpha>\beta\}\cap U\in p$. This complete proof.$\hfill\blacksquare$\\

\begin{exam}
 Let $\{x_\alpha\}$ be a net in topological space $X$. Let $x,y\in X$ be two cluster points of
net. Then there are two subnets $\{x_\beta\}$ and $\{x_\gamma\}$ such that
$x_\beta\rightarrow x$ and $x_\gamma\rightarrow y$. Let $p\in\bigcap_{\beta} cl_{\beta X_d}{\{x_\alpha:\alpha>\beta\}}$ thus by Lemma 3.3
implies that $p-lim x_\alpha=x$ and $p-lim x_\alpha=y$ but $\{x_\alpha\}$ is not convergent.
\end{exam}
\begin{defn}
Let $\{x_\alpha\}_{\alpha\in I}$ be a net in $(X,\tau)$ and $x\in
X$. For $p\in\beta X_d$  we say $Sp-lim x_\alpha=x$ if for each subnet $\{x_\beta\}$ of $\{x_\alpha\}$
implies that $p-lim x_\beta=x$.
\end{defn}

\begin{lem}
Let $S=\{x_\alpha\}_{\alpha\in I}$ be a net in $X$ and $x\in X$.
Then the following statements are equivalent:\\
a) $x_\alpha\rightarrow x$,\\
b) For each $p\in (\bigcap_{\beta\in I}cl_{\beta
X_d}\{x_\alpha:\alpha>\beta\}$, $Sp-lim\,\,x_\alpha=x$.
\end{lem}
{\bf Proof :} $a)\rightarrow b)$. Obvious.\\
$b)\rightarrow a)$. Let $p\in (\bigcap_{\beta\in I}cl_{\beta
X_d}\{x_\alpha:\alpha>\beta\}$, $Sp-lim\,\,x_\alpha=x$. Let $\{x_\alpha\}$ is not convergent to $x$.
Therefore there exists a neighborhood $U$ of $x$ such that for some
subnet $\{x_\beta\}_{\beta \in J}$ of $\{x_\alpha\}$ we have
$\{x_\beta:\beta\in J\}\cap U=\emptyset$. Since $p-lim x_\beta=x$ so for $U\in\tau_x$ and each
$\gamma\in J$, implies that $\{x_\beta:\beta>\gamma \}\cap U\in p$ and this is a contardiction.$\hfill\blacksquare$

\section{\textbf{Generated topology  by ultrafilters}}

Let $X$ be an arbitrary set; by generating a topology on $X$ we mean selecting a family $\tau$ of subsets of $X$
that satisfies conditions of axiom of topology, i.e., a family $\tau$ such that the pair
$(X,\tau)$ is a topological space. We shall now give a new method of generating topology; that consists in
the definition of a neighborhood system combine with ultrafilter converges to a point.
\begin{defn}
For each $x\in X$, let $x^*$ be a subset of $\beta X_d$ such that $\widehat{x}\in x^*$.
For $x\in X$, define $\tau^*_x=\bigcap_{p\in x^*}p$.
\end{defn}

\begin{thm}
a) Suppose we are given a set $X$ and a collection $\{x^*\}_{x\in X}$ of families of subsets of $\beta X_d$. Then
the collection $\{\tau^*_x\}_{x\in X}$ has properties $(BP1)-(BP3)$ in \cite{Engel}. Let $\tau$
be the family of all subsets of $X$ that are unions of subfamilies of $\bigcup_{x\in X}\tau^*_x$. The family
$\tau$ is a topology on $X$ and the collection $\{\tau_x^*\}_{x\in X}$ is a
neighborhood system ( see \cite{Engel}) for the topological
space $(X,\tau)$.\\
The topology $\tau$ is called the topology generated by the ultrafilters system $\{x^*\}_{x\in X}$.\\
b) Let $(X,\tau)$ be a topological space, and $x^*$ be the collection of all ultrafilters
near to $x$ respect to $\tau$. Then the topology generated by the ultrafilters system $\{x^*\}_{x\in X}$
is finer than $\tau$, i.e. $\tau_x\subseteq \tau_x^*$ for each $x\in X$.
\end{thm}
{\bf Proof :} Obvious.$\hfill\blacksquare$
\begin{exam}
 a) For each $x\in X$, let $x^*=\beta X_d$. Then the topology generated by the ultrafilters system
is anti-discrete topology, i.e. $\tau=\{X,\emptyset\}$.

b) For each $x\in X$, let $x^*=\{\widehat{x}\}$, then the topology generated by the ultrafilters system
is discrete topology.

c) Let $X$ be an infinite set and pick $x_\circ\in X$. Define $x^*=\{\widehat{x}\}$ for each $x\neq x_\circ$ and
$x_\circ^*=\beta X_d-\bigcup_{x\neq x_\circ}x^*$. Then all one-point subsets of $X$,
except for the set $\{x_\circ\}$, are open-and-closed; the set $\{x_\circ\}$ is closed but is not open.
For more details see Example 1.1.8 in \cite{Engel}.
\end{exam}
 \begin{thm}[Separation axioms]
 Let $(X,\tau)$ be a topological space and $\{x^*\}_{x\in X}$ be the ultrafilter
 system respect to topology $\tau$. Then :\\
 1) $(X,\tau)$ is a $T_\circ$ space if and only if for each $x,y\in X$, $x^*\neq y^*$.\\
 2) $(X,\tau)$ is a $T_1$ space if and only if for each $x,y\in X$, $x^*- y^*\neq\emptyset$ and $y^*-x^*\neq\emptyset$.\\
 3) $(X,\tau)$ is a $T_2$ space(Hausdorff) if and only if for each $x,y\in X$, $x^*\cap y^*=\emptyset$.\\
 4) $(X,\tau)$ is a regular space if and only if $X$ is a $T_1$ space and for every $x\in X$ and every closed set $F\subseteq X$
 such that $x\notin F$, $x^*\cap cl_{\beta X_d}(F^*)=\emptyset$.\\
 5) $(X,\tau)$ is a completely regular space then for every $x\in X$ and every closed set $F\subseteq X$
 such that $x\notin F$, $cl_{\beta X_d}(x^*)\cap cl_{\beta X_d}(F^*)=\emptyset$.\\
  6)  $(X,\tau)$ is a normal space if and only if $X$ is a $T_1$ space and for every pair
 of disjoint closed subsets $A,B\subseteq X$, $cl_{\beta X_d}(A^*)\cap cl_{\beta X_d}(B^*)=\emptyset$.
 \end{thm}
{\bf Proof :} 1) Let $X$ be a $T_\circ$ space, so for every pair of distinct points $x,y\in X$ there exists
an open set $U$ such that $x\in U$ and $y\notin U$. Thus $x^*\subseteq cl_{\beta X_d}U$ and there exists $q\in y^*$ such that
$q\notin cl_{\beta X_d}U$ and hence $x^*\neq y^*$.

Conversely, let $x^*\neq y^*$ for each pair of distinct points $x,y\in X$. Without loss
of generality, we may assume that $p\in y^*-x^*$, so there is a neighborhood $U$ of $x$ such that $p\notin cl_{\beta X_d}U$.
Notice that $y \notin U$. For, otherwise $y \in U$ implies $p\in y^*\subseteq cl_{\beta X_d}U$ which is a contradiction. So there
exist an open $U$ such that $x\in U$ and $y\notin U$. Thus $X$ is a $T_\circ$ space.

2) Suppose that $X$ is a $T_1$ space and pick $x \neq y$ in $X$. So there exist $U\in \tau_x$ and $V\in\tau_y$ such that
$x\notin V$ and $y\notin U$. Thus $V\notin p$ for some $p\in x^*$ and $U\notin q$ for some $q\in y^*$. Therefore
$x^*-y^*\neq\emptyset$ and $y^*-x^*\neq\emptyset$.

Conversely, let for each $x,y\in X$, $x^*- y^*\neq\emptyset$ and $y^*-x^*\neq\emptyset$. By 1), $x^*- y^*\neq\emptyset$ implies
that there is a neighborhood $U$ of $x$ such that $y\notin U$, and $y^*- x^*\neq\emptyset$ implies
that there is a neighborhood $V$ of $y$ such that $x\notin V$. This complete proof.

3) Let $X$ be a Hausdorff space and pick $x\neq y$ in $X$. So there exist two open sets
$U$ and $V$ such that $x\in U$, $y\in V$ and $U\cap V=\emptyset$. Thus $cl_{\beta X_d}U\cap cl_{\beta X_d}V=\emptyset$. Since
$x^*\subseteq cl_{\beta X_d}U$ and $y^*\subseteq cl_{\beta X_d}V$, so $x^*\cap y^*=\emptyset$.

Conversely, let $X$ not be Hausdorff space, so there exist $x$ and $y$ in $X$ such that $U\cap V\neq\emptyset$ for each
$U\in \tau_x$ and $V\in \tau_y$. Now let
\[
\R_y=\{A\subseteq X:y\in cl_XA\},
\]
so $\R_y$ is partition regular and $\R_y=\R^\uparrow$, see Theorem 2.5. It is obvious that
$\tau_x\subseteq \R_y$. By Theorem 2.5(b), there exists $p\in\beta X_d$ such that $\tau_x\subseteq p\subseteq \R_y$. Therefore
$p\in x^*\cap y^*$ and this is a contradiction.

4) Let $X$ be regular. Let $F$ be closed in $X$ and $x\notin F$, then there are disjoint open sets
$U$ and $V$ in $X$ with $x\in U$ and $F\subseteq V$. Thus
\[
cl_{\beta X_d}U\cap cl_{\beta X_d}V=\emptyset.
\]
This implies that $x^*\cap cl_{\beta X_d}(F^*)=\emptyset$.

Conversely, let $F$ be a closed subset of $X$, $x\notin F$ and let $\tau_F=\{U\in \tau:F\subseteq U\}$. If
$U\cap V\neq\emptyset$  for each $U\in\tau_x$ and $V\in\tau_F$, so $\tau_F\subseteq \R_x$. Thus by Theorem 2.5, there exists
$p\in\beta X_d$ such that $\tau_F\subseteq p\subseteq \R_x$, and hence $p\in x^*\cap \overline{\tau_F}=x^*\cap cl_{\beta X_d}(F^*)=\emptyset$
is a contradiction.

5) Let $X$ be completely regular. Let $F$ be a closed subset of $X$ and pick $x\in X-F$. So there exists a continuous function
$f:X\rightarrow [0,1]$ such that $f|_F=0$ and $f(x)=1$. Hence $f^\beta :B(X)\rightarrow [0,1]$ by $f^\beta(t^*)=f(t)$ for each
$t\in X$ is continuous, see Theorem 3.4. Since $f^\beta$ is continuous , so $(f^\beta)^{-1}(\{0\})$ and $(f^\beta)^{-1}(\{1\})$ are closed,
$F^*\subseteq(f^\beta)^{-1}(\{0\})$ and $(f^\beta)^{-1}(\{0\})\cap (f^\beta)^{-1}(\{1\})=\emptyset$.
Thus $cl_{\beta X_d}(x^*)\cap cl_{\beta X_d}(F^*)=\emptyset$.

6) Let $X$ be normal. Let $A$ and $B$ be two disjoint closed subset of $X$. Then there are disjoint open sets
$U$ and $V$ in $X$ such that $A\subseteq U$ and $B\subseteq V$. Thus
\[
cl_{\beta X_d}(A^*)\cap cl_{\beta X_d}(B^*)\subseteq cl_{\beta X_d}U\cap cl_{\beta X_d}V=\emptyset,
\]
and proof is finish.

Conversely, let $A$ and $B$ be two disjoint closed subset of $X$ such that $cl_{\beta X_d}(A^*)\cap cl_{\beta X_d}(B^*)=\emptyset$.
Define $\tau_A=\{U\in \tau:A\subseteq U\}$ and $\tau_B=\{U\in \tau:B\subseteq U\}$.
If $U\cap V\neq\emptyset$  for each $U\in\tau_A$ and $V\in\tau_B$, so $\tau_B\subseteq \R_A$. Thus by Theorem 2.5, there exists
$p\in\beta X_d$ such that $\tau_B\subseteq p\subseteq \R_A$, and hence $p\in A^*\cap \overline{\tau_B}=A^*\cap cl_{\beta X_d}(B^*)=\emptyset$
is a contradiction. This complete proof.$\hfill\blacksquare$
\begin{defn}
Let $(X,\tau)$ be a topological space. Then \\
a) $X$ is called a $T_{2\frac{1}{2}}$-space if and only if for every $x\neq y$ in $X$,  $x^*\cap \overline{y^*}=\emptyset$.\\
b)  $X$ is called Strongly Hausdorff if and only if for every $x\neq y$ in $X$,  $\overline{x^*}\cap \overline{y^*}=\emptyset$.
 \end{defn}

\begin{exam}
a) Let $X$ be an infinite set, $A$ and $B$ be disjoint subsets of $\beta X_d-X$ such that
$A\cap cl_{\beta X_d}B=\emptyset$ and $cl_{\beta X_d}A\cap cl_{\beta X_d}B\neq\emptyset$. Pick $x_1$ and $x_2$ in $X$, and define
$x_1^*=A\cup\{\widehat{x}_1\}$, $x_2^*=B\cup\{\widehat{x}_2\}$ and $x^*=\{\widehat{x}\}$ for each $x\in X-\{x_1,x_2\}$.
Then $\{x^*\}_{x\in X}$ is an ultrafilter system and let $\tau$ be the topology generated by the ultrafilter system $\{x^*\}_{x\in X}$.
It is obvious that $X$ is a $T_{2\frac{1}{2}}$-space but is not Strongly Hausdorff, because $\{x_2\}$ is closed subset of
$X$ and $x_1\notin \{x_2\}$ and $\overline{x_1^*}\cap \overline{x_2^*}\neq\emptyset$.

b) Let $X$ be an infinite set, $A$ and $B$ be disjoint subsets of $\beta X_d-X$ such that
$cl_{\beta X_d}A\cap cl_{\beta X_d}B\neq\emptyset$. Pick $x_1$ and $x_2$ in $X$, and define
$x_1^*=A\cup\{\widehat{x}_1\}$, $x_2^*=B\cup\{\widehat{x}_2\}$ and $x^*=\{\widehat{x}\}$ for each $x\in X-\{x_1,x_2\}$.
Then $\{x^*\}_{x\in X}$ is an ultrafilter system and let $\tau$ be the topology generated by the ultrafilter system $\{x^*\}_{x\in X}$.
It is obvious that $X$ is Hausdorff. Let $T$ be a closed subset of $X$,
by Lemma 3.1, $x\in T$ if and only if $x^*\cap cl_{\beta X_d}T\neq\emptyset$. Let $x\in X-\{x_1,x_2\}$ and $T$ be a
closed subset of $X$ such that $x\notin T$, then $x^*\cap T^*=\emptyset$.
Now let $x=x_1$ and $T$ be closed subset of $X$ such that $x_1\notin T$ so  $x_1^*\cap T^*=\emptyset$. By similar way, for $x_2\neq x$ and
closed subset $T$ of $X$ that $x_2\notin T$ implies that $x_2^*\cap T^*=\emptyset$. Therefore $X$ is regular space, by Theorem 4.3.
But $X$ is not completely regular and Strongly Hausdorff, because $\{x_2\}$ is closed subset of
$X$, $x_1\notin \{x_2\}$ and $\overline{x_1^*}\cap \overline{x_2^*}\neq\emptyset$.

b) Every Strongly Hausdorff space  is a $T_{2\frac{1}{2}}$ space and hence is Hausdorff.
\end{exam}


\bibliographystyle{alpha}

\end{document}